\documentclass[12pt,a4paper]{article}

\usepackage{amsmath}    % need for subequations
\usepackage{amssymb}
\usepackage{graphicx}   % need for figures
\usepackage{verbatim}   % useful for program listings
\usepackage{color}      % use if colour is used in text

\usepackage{hyperref}   % use for hypertext links, including those to external documents and URLs

\usepackage{enumerate}
\usepackage{mathrsfs}
\usepackage{empheq}
\usepackage{bbold}

\usepackage[final]{showlabels}

\usepackage{amsthm}

\renewcommand{\theta}{\vartheta}
\renewcommand{\phi}{\varphi}

\renewcommand{\title}{Approximate evolution operators for the Active Flux method}

\newcommand{\authorOne}{Wasilij Barsukow\footnote{Max-Planck-Institute for Plasma Physics, Boltzmannstr. 2, 85748 Garching, Germany}}
\begin{document}

\begin{center} \Large
\phantom{m}
\vspace{1cm}

\title

\vspace{1cm}

\date{}
\normalsize

\authorOne
\end{center}

\begin{abstract}

This work focuses on the numerical solution of hyperbolic conservations laws (possibly endowed with a source term) using the Active Flux method. This method is an extension of the finite volume method. Instead of solving a Riemann Problem, the Active Flux method uses actively evolved point values along the cell boundary in order to compute the numerical flux. Early applications of the method were linear equations with an available exact solution operator, and Active Flux was shown to be structure preserving in such cases. For nonlinear PDEs or balance laws, exact evolution operators generally are unavailable. Here, strategies are shown how sufficiently accurate approximate evolution operators can be designed which allow to make Active Flux structure preserving / well-balanced for nonlinear problems.

\end{abstract}

Hyperbolic conservation laws, such as the Euler equations of compressible hydrodynamics, exhibit a large variety of phenomena, in particular in multiple spatial dimensions. This inherent complexity is reflected in high computational cost associated with attempts of solving the equations numerically. However, often additional difficulties arise because of inefficient or unadapted numerical methods, i.e. methods which require a resolution (in time or space) that is much higher than the desired resolution of the sought numerical solution. It is this latter aspect that is addressed in this work, which focuses on the development of a novel family of numerical methods (\emph{Active Flux}).

With the Cauchy problem being the natural setting for hyperbolic conservation laws, it is customary to integrate the numerical approximations forward in time in an explicit way. This, however, turns out to require \emph{upwinding}, which replaces (unstable) central spatial derivatives by specific one-sided derivatives. The choice of upwind direction traditionally (both for Finite Volume and Finite Element methods, such as Discontinuous Galerkin) involves Riemann solvers, i.e. exact or approximate solutions of initial value problems associated to discontinuous data. Godunov's method proposes to introduce such a discontinuity at every cell interface. 

The extension of this strategy to multiple spatial dimensions does not reach up to the success of some finite difference methods (e.g. \cite{morton01,barsukow17a}) that were obtained by a less fundamental approach, as was shown in \cite{barsukow17}. These methods (examples of so-called \emph{structure preserving} methods) preserve discrete involutions and discretize all the stationary states of a PDE, instead of keeping stationary discretizations of merely trivial ones. All this is associated with the practical advantage of achieving excellent results on coarse grids. It thus has become clear that there is a large potential of saving computational cost, but a lack of fundamental principles according to which such methods can be derived.

\section{The Active Flux method}

A promising structure preserving numerical method is the novel Active Flux method (\cite{vanleer77,eymann13}), which presents an alternative to the traditional Godunov idea. It has been initially derived for linear problems, and was shown to be vorticity and stationarity preserving for linear acoustics on Cartesian grids (\cite{barsukow18activeflux}). The idea of Active Flux is to reconstruct the solution globally continuously, by introducing additional pointwise degrees of freedom at cell boundaries. It thus is a blend of Finite Volume and Finite Element ideas. The evolution of the cell averages follows trivially from the available data at cell boundaries. The evolution of the pointwise degrees of freedom, however, does not follow the usual Finite Element approach, which might be in parts the reason behind Active Flux' success.

Finite element (and also finite volume) methods often adapt a method-of-lines approach, separating the spatial discretization from the integration in time. Performing the former, one ends up with a large system of ODEs, which are solved by more or less standard ODE integrators, yielding the latter. On the contrary, in the initial versions of Active Flux, devoted to linear problems, the point values were updated by applying the exact evolution operator to the initial data given by the (piecewise parabolic) continuous reconstruction of the data at the previous time step. 

This would mean that the point value $q_{i+\frac12}^{n+1}$ at time $t^{n+1}$ is the value $q(t^{n+1} - t^n,x_{i+\frac12} )$ of an exact evolution $q(t,x)$ of piecewise defined initial data
\begin{align}
q(0, x) &= -3 (2 \bar q_i^n - q^n_{i-\frac12} - q^n_{i+\frac12}) \frac{(x-x_i)^2}{\Delta x^2} \\&+ (q^n_{i+\frac12} - q^n_{i-\frac12}) \frac{x-x_i}{\Delta x} + \frac{6 \bar q_i^n - q^n_{i-\frac12} - q^n_{i+\frac12}}{4}  \quad \text{if} \quad x \in [x_{i-\frac12},x_{i+\frac12}]
\end{align}
where $\bar q_i^n$ denotes the average of $q$ in cell $i$ at time $t^n$. Note that
\begin{align}
 q(0,x_{i\pm\frac12}) &= q^n_{i\pm\frac12} & \frac{1}{\Delta x} \int_{x_{i-\frac12}}^{x_{i+\frac12}} q(0,x) \,\mathrm{d}x &= \bar q_i^n
\end{align}

For linear advection $\partial_t q + c \partial_x q = 0$, the evolution operator is just $q(t,x) = q(0, x- ct)$, and thus (with the CFL number $\lambda = \frac{c(t^{n+1}-t^n)}{\Delta x}$) one obtains for $c> 0$
\begin{align}
q_{i+\frac12}^{n+1} &= -6 \bar q_i^n ( \lambda-1 )\lambda + q_{i-\frac12}^n\lambda ( 3 \lambda-2 )  + q_{i+\frac12}^n (\lambda-1)( 3 \lambda-1 ) \label{eq:timeevolin}\\
\bar q^{n+1}_i &= \bar q_i^n - \lambda \frac{q_{i+\frac12}^n + 4 q_{i+\frac12}^{n+\frac12} + q_{i+\frac12}^{n+1} - q_{i-\frac12}^n + 4 q_{i-\frac12}^{n+\frac12} + q_{i-\frac12}^{n+1}}{6}
\end{align}
which cannot be obviously associated to a spatial discretization and a subsequent application of an ODE solver (compare to examples in \cite{abgrall20}). 

The advantage of time evolution \eqref{eq:timeevolin} lies in its immediate stability, ranging up to the physical stability condition of $\lambda \leq 1$, and a natural derivation. Interestingly, it has been shown in \cite{kerkmann18} that the above method can be given the interpretation of an ADER-type evolution: the formal Taylor series of $q(t, x)$ in time at $x=x_{i+\frac12}$ is truncated and the time derivatives replaced by spatial derivatives using the PDE. Whereas the value of $q$ at $x_{i+\frac12}$ does exist, the globally continuous reconstruction does not imply the same for derivatives. In \cite{kerkmann18} it is shown that for linear one-dimensional problems replacing the evaluation of a derivative at the location of its discontinuity by solutions of Riemann problems \emph{in the derivatives} yields exactly the same method \eqref{eq:timeevolin}. However, for nonlinear problems the properties of such an approach are less clear, the equation governing the evolution of derivatives complicated, and the solution of Riemann problems for these (non-conservative) equations, generally unknown.

The work presented here is part of an effort to extend Active Flux to nonlinear problems in a way that stays as close possible to the approach used for linear problems. More precisely, the aim is to find (iterative) approximations to evolution operators, such that they revert to exact evolution operators on linear problems. For example, in case of a scalar nonlinear conservation law, the concept of characteristics persists from the linear case, and hope is that the speed of the characteristic can be estimated to sufficient accuracy. Recent advances for this case, and the case of one-dimensional systems, are presented next.

\section{Scalar conservation laws}

An approximate evolution operator for scalar conservation laws has been suggested in \cite{barsukow19activeflux}. Assume for the moment that no shocks are present. For 
\begin{align}
 \partial_t q + a(q) \partial_x q &= 0 & a&\colon \mathbb R \to \mathbb R & q&\colon \mathbb R^+_0 \times \mathbb R \to \mathbb R
\end{align}
the characteristic starting at $x = \xi$ at time $t=0$ is subject to the relation
\begin{align}
 \xi = x - a(q(0,\xi)) t
\end{align}
which can be solved iteratively:
\begin{align}
 \xi^{(0)} &:= x  & \xi^{(n+1)} &:= x - a(q(0,\xi^{(n)})) t \label{eq:timeevononlin}
\end{align}
$\xi^{(n)}$ is approximating $\xi$ to an error $\mathcal O(t^{n+1})$ (\cite{barsukow19activeflux}), and moreover for linear problems, one iteration gives the exact solution.

\section{Systems of conservation laws in 1-d}

For the case of systems in one spatial dimension the concept of characteristics still persists, but they are curved and, generally speaking, no quantities are constant along any of them. For simplicity, however, assume first that such quantities exist (as is the case for the shallow water equations, for example). Then, a nonlinear $m \times m$ system
\begin{align}
 \partial_t q + J(q) \partial_x q &= 0 & q \colon \mathbb R^+_0 \times \mathbb R \to \mathbb R^m\\
\intertext{can be written as}
 \partial_t Q + \mathrm{diag}(\lambda_1, \ldots, \lambda_m) \partial_x Q &= 0 & Q \colon \mathbb R^+_0 \times \mathbb R \to \mathbb R^m
\end{align}
where each eigenvalue $\lambda_i$ of $J$ is considered a function of $Q = (Q_1, \ldots, Q_m)^\text T$. Thus, $Q_i$ is constant along the $i$-th characteristic, but the characteristic is curved, as its slope additionally depends on the values of the other variables.

Considering the first two steps of the fixpoint iteration that successfully yielded an approximate evolution operator in the scalar case, one might be tempted to choose ($m=2$ to save a tree)
\begin{align}
 \tilde Q_j(t, x) = Q_{j,0}\Big(x - t\lambda_j(Q_{1,0}(x - t\lambda_j(x)), Q_{2,0}(x - t \lambda_j(x))) \Big) \quad j=1,2
\end{align}
as an approximate evolution operator. Here, $Q_{j,0}$ denotes the initial data of $Q_j$ and $\lambda_j(x)$ is short-hand for $\lambda_j(Q_{1,0}(x), Q_{2,0}(x))$. Unfortunately, the error (in $t$) of this approximation is no better than simply taking
\begin{align}
 \tilde Q_j(t, x) = Q_{j,0}\Big(x - t\lambda_j(x) \Big)
\end{align}
and not enough to achieve a third order method, as is customary for Active Flux.

This possibly surprising result has to do with the fact that characteristics are curved. It turns out that an approximate operator of sufficient order of accuracy is obtained as follows:
\begin{align}
 \tilde Q_j(t, x) &= Q_{j,0}\left(x - t\lambda_j\left(Q_{1,0}\left(x - t\frac{\lambda_1(x) + \lambda_j(x)}{2}\right),\right.\right.\\
 &\phantom{mmmmmmmmll}\left.\left.Q_{2,0}\left(x - t \frac{\lambda_2(x) + \lambda_j(x)}{2}\right)\right) \right) \quad j=1,2
\end{align}
Note how it reduces to the second iteration of operator \eqref{eq:timeevononlin} in the scalar case. 

The reason for this is the following. Imposing the evolution operator to have the shape $\tilde Q_j(t, x) = Q_{j,0}\Big(x - t\lambda_j^* \Big)$ for some $\lambda_j^*$ means that the curved characteristic is replaced by a straight characteristic with \emph{average} speed $\lambda^*$, which in general is different both from its speed at the footpoint and at time $t$. Thus,
\begin{align}
 \lambda_j^* \simeq \frac{1}{t} \int_0^t  \lambda_j\Big(Q_1(\tau, X_j(\tau)), \ldots , Q_m(\tau, X_j(\tau))\Big) \mathrm{d} \tau
\end{align}
where $X_j\colon \mathbb R^+_0 \to \mathbb R$ denotes the $j$-th characteristic curve. It is not surprising (and is proved to be true in \cite{barsukow19activeflux}) that a good estimate therefore is
\begin{align}
 \lambda_j^* &\simeq \lambda_j \left(Q_1\left(\frac{t}{2}, X_j\left(\frac{t}{2}\right)\right), \ldots , Q_m\left(\frac{t}{2}, X_j\left(\frac{t}{2}\right)\right) \right) 
\end{align}
and 
\begin{align}
 Q_i\Big(\frac{t}{2}, X_j\Big(\frac{t}{2}\Big)\Big) \simeq  Q_{i}\Big(\frac{t}{2}, x - \frac{t}{2} \lambda_j(x)\Big) \simeq Q_{i,0}\Big(x - t \frac{\lambda_i(x) + \lambda_j(x)}{2} \Big)
\end{align}
which is the expression stated earlier.

For systems without characteristic variables (such as the Euler equations) simi\-lar approximate evolution operators can be found. In that case, not only the characteristic speeds need to be approximated, but also the transformation matrix that diagonalizes $J$. For details, see \cite{barsukow19activeflux}.

Both for scalar problems and for systems it is necessary to anticipate the case of crossing characteristics, and thus shock formation. Because of global continuity of the reconstruction, a shock does not form instantaneously, and often enough, the CFL condition imposes a smaller time step than the shock formation time. In other cases, it is suggested in \cite{barsukow19activeflux} to compute several estimates of possible characteristics, and select one of them.

\section{Outlook}

Future work will be devoted to the multi-dimensional case. There, the additional difficulty is the conceptual replacement of characteristic lines by characteris\-tic cones. This means that predictor-corrector strategies are not only needed to approximate the speed (i.e. the tangent to the cone), but also the directions.


\begin{thebibliography}{BHKR19}

\bibitem[Abg20]{abgrall20}
R{\'e}mi Abgrall.
\newblock A combination of residual distribution and the active flux
  formulations or a new class of schemes that can combine several writings of
  the same hyperbolic problem: application to the 1d euler equations.
\newblock {\em arXiv preprint arXiv:2011.12572}, 2020.

\bibitem[Bar19]{barsukow17a}
Wasilij Barsukow.
\newblock Stationarity preserving schemes for multi-dimensional linear systems.
\newblock {\em Mathematics of Computation}, 88(318):1621--1645, 2019.

\bibitem[Bar21]{barsukow19activeflux}
Wasilij Barsukow.
\newblock The active flux scheme for nonlinear problems.
\newblock {\em Journal of Scientific Computing}, 86(1):1--34, 2021.

\bibitem[BHKR19]{barsukow18activeflux}
Wasilij Barsukow, Jonathan Hohm, Christian Klingenberg, and Philip~L Roe.
\newblock The active flux scheme on {C}artesian grids and its low {M}ach number
  limit.
\newblock {\em Journal of Scientific Computing}, 81(1):594--622, 2019.

\bibitem[BK20]{barsukow17}
Wasilij Barsukow and Christian Klingenberg.
\newblock Exact solution and a truly multidimensional {G}odunov scheme for the
  acoustic equations.
\newblock {\em submitted, preprint available as arXiv:2004.04217}, 2020.

\bibitem[ER13]{eymann13}
Timothy~A Eymann and Philip~L Roe.
\newblock Multidimensional active flux schemes.
\newblock In {\em 21st AIAA computational fluid dynamics conference}, 2013.

\bibitem[HKS19]{kerkmann18}
Christiane Helzel, David Kerkmann, and Leonardo Scandurra.
\newblock A new {ADER} method inspired by the active flux method.
\newblock {\em Journal of Scientific Computing}, 80(3):1463--1497, 2019.

\bibitem[MR01]{morton01}
Keith~William Morton and Philip~L Roe.
\newblock Vorticity-preserving {L}ax-{W}endroff-type schemes for the system
  wave equation.
\newblock {\em SIAM Journal on Scientific Computing}, 23(1):170--192, 2001.

\bibitem[vL77]{vanleer77}
Bram van Leer.
\newblock Towards the ultimate conservative difference scheme. {IV}. {A} new
  approach to numerical convection.
\newblock {\em Journal of computational physics}, 23(3):276--299, 1977.

\end{thebibliography}
\end{document}